\documentclass[12pt]{article}
\usepackage{amsmath,amsfonts,amssymb,rotating,colordvi}
\usepackage{color}
\usepackage[unicode=true,
 bookmarks=true,bookmarksnumbered=true,bookmarksopen=true,bookmarksopenlevel=2,
 breaklinks=false,pdfborder={0 0 1},backref=false,colorlinks=true]
 {hyperref}

\usepackage{tikz}
\usetikzlibrary{arrows,shapes,positioning}
\usetikzlibrary{decorations.markings}
\tikzstyle arrowstyle=[scale=1]
\tikzstyle directed=[postaction={decorate,decoration={markings,mark=at position .65 with {\arrow[arrowstyle]{stealth}}}}]
\tikzstyle reverse directed=[postaction={decorate,decoration={markings,mark=at position .65 with {\arrowreversed[arrowstyle]{stealth};}}}]

\hypersetup{pdftitle={The unimodality of the Ehrhart
$\delta$-polynomial of
the chain polytope of the zig-zag poset},
 pdfauthor={ },
 pdfsubject={ },
 pdfkeywords={ },
 linkcolor=black,citecolor=black,urlcolor=blue,filecolor=blue,pdfpagelayout=OneColumn,pdfnewwindow=true,pdfstartview=XYZ,plainpages=false}

\def\qed{\nopagebreak\hfill{\rule{4pt}{7pt}}}

 \newtheorem{thm}{Theorem}[section]

\newtheorem{conj}[thm]{Conjecture}

\numberwithin{equation}{section}

\newdimen\Squaresize \Squaresize=11pt
\newdimen\Thickness \Thickness=0.7pt
\def\Square#1{\hbox{\vrule width \Thickness
   \vbox to \Squaresize{\hrule height \Thickness\vss
    \hbox to \Squaresize{\hss#1\hss}
   \vss\hrule height\Thickness}
\unskip\vrule width \Thickness} \kern-\Thickness}

\def\Vsquare#1{\vbox{\Square{$#1$}}\kern-\Thickness}

\def\moins{\raise 1pt\hbox{{$\scriptstyle -$}}}

\begin{document}

\begin{center}
\textbf{\large{}The unimodality of the Ehrhart
$\delta$-polynomial\\
 of
the chain polytope of the zig-zag poset}\textbf{ }
\par\end{center}

\begin{center}
Herman Z.Q. Chen$^{1}$ and Philip B. Zhang$^{2}$\\[6pt]
\par\end{center}

\begin{center}
$^{1}$Center for Combinatorics, LPMC\\
 Nankai University, Tianjin 300071, P. R. China

$^{2}$College of Mathematical Science \\
Tianjin Normal University, Tianjin  300387, P. R. China\\[6pt]
\par\end{center}

\begin{center}
Email: $^{1}$\texttt{zqchern@163.com},
$^{2}$\texttt{zhangbiaonk@163.com}
\par
\end{center}

\noindent \textbf{Abstract.}
We prove the unimodality of the Ehrhart $\delta$-polynomial of the chain polytope of the
zig-zag poset, which was conjectured by Kirillov. First, based on a result due
to Stanley, we show that this polynomial coincides with the $W$-polynomial for
the zig-zag poset with some natural labeling. Then, its unimodality immediately
follows from a result of Gasharov, which states that the $W$-polynomials of
naturally labeled graded posets of rank $1$ or $2$ are unimodal.

\noindent \emph{AMS Classification 2010:} 05A15, 05A20


\section{Introduction}

The main objective of this paper is to prove a unimodality conjecture on
$\delta$-polynomials of the chain polytope of the zig-zag poset,
which was proposed by Kirillov \cite{Kirillov99} in the study of Kostka numbers and Catalan
numbers. Let us first give an overview of Kirillov's conjecture.

Let $\mathbb{Z}^m$ denote the $m$-dimensional integer lattice
in $\mathbb{R}^m$, and let $\mathcal{P}$
be an $m$-dimensional lattice polytope in $\mathbb{R}^m$.
A remarkable theorem due to Ehrhart \cite{Ehrhart62} states that
the number of
lattice points that lie inside the dilated polytope $n\mathcal{P}$:
\begin{align}
i(\mathcal{P};n)= |n\mathcal{P}\cap \mathbb{Z}^m| .
\end{align}
is given by a polynomial in $n$ of
degree $m$, called the Ehrhart polynomial of the lattice polytope
$\mathcal{P}$. By a well known result about rational generating functions, see
\cite[Corollary 4.3.1]{Stanley96}, the generating
function (called {the Ehrhart series of $\mathcal{P}$})
\begin{equation}\label{J_i_delta}
J(\mathcal{P};t)=\sum_{n\geq 0}i(\mathcal{P};n)t^n
\end{equation}
evaluates to a rational function:
\begin{equation}\label{J_i_delta}
J(\mathcal{P};t)=\frac{\delta(\mathcal{P}
;t)}{(1-t)^{\dim{\mathcal{P}+1}}}
\end{equation}
for some polynomial $\delta(\mathcal{P};t)$ of degree at most
$\dim(\mathcal{P})$, which is called
{the Ehrhart $\delta$-polynomial of $\mathcal{P}$}.
If the polynomial $\delta(\mathcal{P};t)$ is of the following form
$$\delta(\mathcal{P};t)=\delta_0+\delta_1 x+\cdots+\delta_m x^m,$$
then we call $(\delta_0, \delta_1, \ldots, \delta_m)$ the {(Ehrhart) $\delta$-vector} of $\mathcal{P}$.
Stanley \cite{Stanley80} also proved that $\delta(\mathcal{P};t)$ must be a polynomial in nonnegative coefficients. For more information on the Ehrhart theory of rational polytopes, see \cite{Beck07}.

Let $\mathcal{P}_n$ be a convex integral polytope in $\mathbb{R}^n$ determinated by the following
inequalities
\[
\begin{array}{cl}
x_i\geq 0,		& \mbox{for}\, 1\leq i\leq n, \\
x_i+x_{i+1}\leq1,	& \mbox{for}\, 1\leq i\leq n-1.
\end{array}
\]
Kirillov conjectured that the $\delta$-polynomial 
of $\mathcal{P}_n$ is unimodal \cite{Kirillov99}. Recall that a polynomial $f(x)=\sum_{i=0}^{n}a_i x^i$ with real coefficients is
said to be {unimodal} if there exists an integer $i\geq0$
such that
$$a_0\leq \cdots \leq a_{i-1}\leq a_i \geq a_{i+1}\geq\cdots\geq a_n,$$
and {symmetric} if
for all $0\leq i \leq n$
$$a_i=a_{n-i}.$$ Kirillov's conjecture is stated as follows.

\begin{conj}[{\cite[p.119, Conjecture 3.11]{Kirillov99}}]\label{kirillov-conj}
For any $n\geq 1$, the $\delta$-polynomial $\delta(\mathcal{P}_n;t)$ is unimodal.
\end{conj}

In this paper, we give a proof of Kirillov's conjecture. Our proof is based on the theory of chain polytopes of posets, as well as the theory of $W$-polynomials of posets.

\section{Preliminaries}

In this section, we shall review some definitions and results on chain polytopes and $W$-polynomials of posets.

We begin with some definitions concerning posets. Let $(P,\preceq)$ be a poset with $d$ elements. Recall that
a {chain} of length $\ell$ in $P$ is a sequence $a_0\prec a_1\prec\cdots\prec
a_{\ell}$, and it is called {maximal} in $P$ if we cannot add elements to this chain.
If every maximal chain of $P$ has the same length $r$,
then we say that $P$ is {graded} of {rank} $r$ and denote the rank of $P$ by $\mathrm{rank}(P)$.
In this case, there is a unique {rank function} $\rho: P\rightarrow \{0,1,\ldots,d\}$
such that $\rho(x)=0$ if $x$ is a minimal element of $P$, and $\rho(y)=\rho(x)+1$ if $y$ covers
$x$ in $P$. If $\rho(x)=i$, then we say that $x$ is of rank $i$.

The notion of chain polytopes was introduced by Stanley \cite{Stanley86}. Given a poset $P$ with elements $\{a_1,\ldots,a_d\}$, Stanley associated it with a polytope $\mathcal{C}(P)$ defined by the chains in $P$, called the {chain polytope} of $P$. Precisely, the chain polytope $\mathcal{C}(P)$ is the convex polytope consisting of those
$(x_1,\ldots,x_d)\in \mathbb{R}^d$ such that
\begin{itemize}
 \item $x_i\geq 0$, {for every} $a_i\in P$,
 \item $x_{p_1}+x_{p_2}+\cdots+x_{p_k}\leq1$, {for every
chain} $a_{p_1}\prec \cdots\prec a_{p_k}$ of $P$.
\end{itemize}
Since $\mathcal{C}(P)$ contains the $d$-dimensional simplex
$$\{(x_1,\ldots,x_d)\in \mathbb{R}^d : x_i\geq 0 \mbox{ for all } 1\leq i \leq
d \mbox{ and } x_1+x_2+\cdots+ x_d \leq 1\},$$ 
we know that $$\dim (\mathcal{C}(P))=d=|P|.$$

We would like to point out that the polytope $\mathcal{P}_n$ is just the chain
polytope of the zig-zag
poset of order $n$. Recall that the {zig-zag
poset} of order $n$ is the poset
$\mathbf{Z}_n=\{a_1,\ldots,a_n\}$ in which 
$$a_1\prec a_2\succ a_3\prec\cdots,$$ 
see \cite{Blyth94,Stanley96}.
Note that, if $n=2k+1$, the maximal chains of $\mathbf{Z}_n$ are
$a_{2i-1}\prec a_{2i}$ and $a_{2i+1}\prec a_{2i}$ for $1\leq i\leq k$. While, if $n=2k+2$,
the maximal chains of $\mathbf{Z}_n$ are $a_{2i-1}\prec a_{2i}$ and $a_{2i+1}\prec a_{2i}$ for $1\leq i\leq k$, together with $a_{2k+1}\prec a_{2k+2}$. By definition, it is clear that
\begin{align}\label{p=z}
\mathcal{P}_n=\mathcal{C}(\mathbf{Z}_n).
\end{align}

Based on the above viewpoint, Conjecture \ref{kirillov-conj} is equivalent to the statement that the
$\delta$-polynomial of $\mathcal{C}(\mathbf{Z}_n)$ is unimodal. While for the chain polytope of poset $P$,
Stanley \cite{Stanley86} has already established a connection between the $\delta$-polynomial of the chain polytope and the number of order-preserving maps of $P$. Let $m$
be a positive integer and let
$\widetilde{\Omega}(P;m)$ denote the number of
order-preserving maps $\eta:P\rightarrow\{1,2,\ldots,m\}$,
i.e., if $x\preceq y$ in $P$ then $\eta(x)\leq \eta(y)$.
It is known that
$\widetilde{\Omega}(P;m)$ is a polynomial
of degree $|P|$ in $m$. Equivalently, there exists a polynomial
$\widetilde{W}(P;t)$ of degree $\leq
|P|$ such that
\begin{align}\label{Omega_WP}
 \sum_{m\geq 0}
\widetilde{\Omega}(P;m+1)t^m=\frac{\widetilde{W}(P;t)}{(1-t)^{|P|+1}}.
\end{align}
Stanley obtained the following theorem.

\begin{thm}[{\cite[Theorem 4.1]{Stanley86}}]\label{i_Omega}
For any positive integer $m$ and any poset $P$, we have
 $$i(\mathcal{C}(P);m)=\widetilde{\Omega}(P;m+1),$$
 or equivalently,
 \begin{align}\label{cp=pre}
 \delta(\mathcal{C}(P);t)=\widetilde{W}(P;t).
 \end{align}
\end{thm}

Instead of considering the number of order-preserving maps of $P$, we may also
study the number of order-reversing maps of $P$. In fact, there is a more general theory on
order-reversing maps, developed by Stanley \cite{Stanley72} and called the theory of $P$-partitions.
Suppose that $P$ is a
finite poset
with $d$ elements as before. A labeling $\omega$ of $P$ is a bijection from $P$ to
$\{1,2,\ldots, d\}$. The labeling $\omega$ is called natural if $x\preceq y$
implies $\omega(x)\leq \omega(y)$ for any $x,y\in P$, namely, it is an order-preserving map.
A {$(P, \omega)$-partition} is a map $\sigma$ which
satisfies the following conditions:

\begin{itemize}
\item $\sigma$ is order reversing, namely, $\sigma(x)\geq \sigma(y)$ if
$x\preceq y$ in $P$; and moreover
\item if $\omega(x)>\omega(y)$, then $\sigma(x)>\sigma(y)$.
\end{itemize}
The {order polynomial} $\Omega(P,\omega; n)$ is defined as the number of $(P,
\omega)$-partitions
$\sigma$ with $\sigma(x)\leq n$ for any $x\in P$. It is also known that
$\Omega(P,\omega; n)$ is a polynomial
of degree $|P|$ in $n$, or equivalently, there exists a polynomial
$W(P,\omega;t)$, called the $W$-polynomial of $(P,\omega)$, of degree $\leq
|P|$ such that
\begin{align}\label{Omega_WP}
 \sum_{n\geq 0} \Omega(P,\omega;n+1)t^n=\frac{W(P,\omega;t)}{(1-t)^{|P|+1}}.
\end{align}

Note that, for a natural labeling $\omega$, we must have
\begin{align}\label{rev=pre}
\Omega(P,\omega;n)=\widetilde{\Omega}(P;n),
\end{align}
since $\Omega(P,\omega;n)$ is just the number of order-reversing maps in this case,
and $\widetilde{\Omega}(P;n)$ is the number of order-preserving maps.
In fact, if $\eta:P\rightarrow\{1,2,\ldots,m\}$ is order reversing, then 
the map $\widetilde{\eta}:P\rightarrow\{1,2,\ldots,m\}$ defined by
$$\widetilde{\eta}(x)=m+1-\eta(x)$$
is order-preserving, and vice versa. 
By \eqref{cp=pre} and \eqref{rev=pre}, we have
\begin{align}\label{key-id}
\delta(\mathcal{C}(P);t)=\widetilde{W}(P;t)=W(P,\omega;t).
\end{align}

\section{Proof }

In this section, we shall give a proof of Conjecture \ref{kirillov-conj}.
Our proof is based on the following result due to Gasharov \cite{Gasharov98}.

\begin{thm}[{\cite[Theorem 1.2]{Gasharov98}}]\label{w-unimodal}
 If $P$ is a graded poset with $1\leq \mathrm{rank}(P) \leq 2$ and $\omega$ is a
natural labeling of
 $P$, then $W(P,\omega;t)$ is unimodal.
\end{thm}

We proceed to prove Conjecture \ref{kirillov-conj}.

\noindent\textit{Proof of Conjecture \ref{kirillov-conj}.}
By \eqref{key-id}, we have
$$\delta(\mathcal{C}(\mathbf{Z}_n);t)=\widetilde{W}(\mathbf{Z}_n;t)=W(\mathbf{Z}_n,\omega;t)$$
for some natural labeling $\omega$ of the zig-zag poset $\mathbf{Z}_n$. It is clear that
$\mathbf{Z}_n$ a graded poset with $\mathrm{rank}(P)=1$. From Theorem
\ref{w-unimodal} it follows the
unimodality of $\delta(\mathcal{C}(\mathbf{Z}_n);t)$, and hence that of $\delta(\mathcal{P}_n;t)$.
This completes the proof.
\qed

\vskip 3mm
\noindent{\bf Acknowledgements.} This work was supported by the 973 Project, the PCSIRT Project of the Ministry of Education and the National Science Foundation of China.

\end{document}